\documentclass[11pt,reqno]{amsart}

\usepackage[T1]{fontenc}
\usepackage{lmodern}
\usepackage{amssymb}
\usepackage{microtype}
\usepackage[a4paper,margin=29mm]{geometry}
\usepackage{xcolor}
\usepackage[colorlinks,allcolors=blue!55!black]{hyperref}

\newtheorem{theorem}{Theorem}[section]
\newtheorem{proposition}[theorem]{Proposition}
\newtheorem{lemma}[theorem]{Lemma}
\DeclareMathOperator{\lcm}{lcm}
\newcommand{\Z}{\mathbb Z}
\newcommand{\R}{\mathbb R}
\newcommand{\T}{\mathbb T}
\newcommand{\E}{\mathbb E}
\newcommand{\ind}{\mathbf 1}
\numberwithin{equation}{section}
\allowdisplaybreaks[2]
\title{A logarithmic saving in the sum-free subset problem}
\author{Xiaolu Zhang}
\email{zxl0714@gmail.com}
\hypersetup{pdftitle={A logarithmic saving in the sum-free subset problem},
  pdfauthor={Xiaolu Zhang},pdfsubject={Erdos Problem 792 and sum-free subsets}}

\begin{document}

\begin{abstract}
For every sufficiently large \(n\) there is a set of \(n\) positive
integers, all at most \(n^2\), whose largest sum-free subset has size
at most \(n/3+O(n/\log n)\). We bound the number of small divisible
elements in a dense set with few Schur triples. A fixed-accuracy
cyclic model gives a rotation, and the original integer triple count
controls its density at the required scale.
\end{abstract}

\maketitle

\section{Introduction}\label{sec:intro}

A set of integers is \emph{sum-free} if it contains no solution to
\(x+y=z\), allowing \(x=y\). For a finite integer set \(A\), let
\(\alpha(A)\) be the largest size of a sum-free subset of \(A\), and put
\[
 f(n)=\min_{\substack{A\subseteq\Z\\|A|=n}}\alpha(A).
\]
Erd\H{o}s Problem~\#792 asks for the order of \(f(n)-n/3\) up to
constant factors \cite{bloom}. Erd\H{o}s proved the lower bound
\(n/3\) for sets of \(n\) nonzero integers \cite{erdos}.
Eberhard, Green and Manners proved \(f(n)\le n/3+o(n)\)
using positive integers \cite[Theorem~1.1]{egm}. Their bound also
holds for subsets that avoid only triples with three distinct entries.
Eberhard gave a finitary approach using multiplicative F{\o}lner
sets \cite{eberhard}.

The question of an unbounded lower excess, also recorded by Green
\cite[Problem~1]{green}, was answered by Bedert. His result gives
\(f(n)\ge n/3+c\log\log n\) for some \(c>0\) and all sufficiently
large \(n\) \cite[Theorem~1.2]{bedert}. Deleting zero, if present,
changes the lower bound by only \(O(1)\).

\begin{theorem}\label{thm:main}
There are absolute constants \(C\) and \(n_0\) such that for each
integer \(n\ge n_0\) there is a set \(A\subseteq[1,n^2]\cap\Z\) with
\[
 |A|=n,\qquad \alpha(A)\le \frac n3+C\frac{n}{\log n}.
\]
\end{theorem}

In particular \(f(n)-n/3=O(n/\log n)\). The order of the excess between
this upper bound and the lower bound of Bedert remains open.

\subsection*{Outline of the proof}
We follow the weighting strategy of Eberhard, Green and Manners
\cite[Sections~3 and~5]{egm}. The extra weight goes on small multiples
of a modulus, where a dense set with few Schur triples has little
mass. Proposition~\ref{prop:shifts} gives an interval of length
proportional to \(\eta^2 N\) and a modulus of size \(\exp(O(1/\eta))\)
at density surplus \(\eta\). Iterating at geometrically decreasing
surpluses gives a logarithmic cost \(O(1/\epsilon)\) to reach error
\(\epsilon\). Sampling at cardinality \(n\) then allows
\(\epsilon\) of order \(1/\log n\).

Removing the lower endpoints of pairs at a prescribed distance
gives a set disjoint from its translate. Sections~\ref{sec:sampling}
and~\ref{sec:structure} approximate this set at fixed accuracy by
interval slices on a cyclic quotient, with centers moving linearly
in time. Sections~\ref{sec:return-integers} and~\ref{sec:integer-density}
use the original ordered triple count to control the phase offset
of the resulting rotation and the density near the two ends of
the middle third. These integer estimates give the dependence on
\(\eta\). The stability theorems enter at a fixed accuracy.

\section{From small shifts to finite witnesses}\label{sec:witnesses}

Write \([N]=\{1,\dots,N\}\). For \(A\subseteq[N]\), define
\[
 T(A)=|\{(x,y)\in A^2:x+y\in A\}|,\qquad I_N(A)=\frac{T(A)}{N^2}.
\]
The count is ordered and includes equal inputs. Indicators of
integer sets are extended by zero to \(\Z\). We use \(\E\) for
averages with respect to uniform measures on finite sets and Haar
probability measure on \(\T=\R/\Z\).

For a positive shift \(t\), let \(L_A(t)=|\{x\in A:x+t\in A\}|\).
If \(t\in A\), each pair counted by \(L_A(t)\) gives the ordered
triple \((t,x,t+x)\). Summing over its first input gives
\begin{equation}\label{eq:shift-triples}
 T(A)=\sum_{t\in A}L_A(t).
\end{equation}

\begin{proposition}\label{prop:shifts}
There are absolute constants \(K,a,\eta_0>0\), with \(\eta_0\le1/12\),
and an integer \(H_0\ge1\) such that the following holds. For
\(0<\eta\le\eta_0\), set
\[
 H_\eta=\lceil H_0/\eta\rceil,\qquad Q_\eta=\lcm(1,\dots,H_\eta),
 \qquad \kappa_\eta=K\eta^2\le \eta/16.
\]
For every \(N\), if \(A\subseteq[N]\), \(|A|\ge(1/3+\eta)N\), and
\(I_N(A)\le a\eta^2\), then
\begin{equation}\label{eq:origin}
 |A\cap Q_\eta\Z\cap[1,\kappa_\eta N]|\le \frac{8 T(A)}{\eta N}.
\end{equation}
Moreover \(\log Q_\eta=O(1/\eta)\).
\end{proposition}

Under the hypotheses of the proposition, Section~\ref{sec:proof-shifts}
shows that \(L_A(t)>\eta N/8\) for every positive multiple \(t\)
of \(Q_\eta\) up to \(\kappa_\eta N\). Each such multiple that lies in \(A\)
therefore accounts for more than \(\eta N/8\) triples in
\eqref{eq:shift-triples}. A dense set with few triples contains
few points of this progression, which is where the construction
puts its extra weight.

\subsection{Weights on integer intervals}

\begin{lemma}\label{lem:weights}
For each \(0<\epsilon\le1/12\) there are \(N_*,W,\delta^{-1}
\le \exp(C_0/\epsilon)\), with \(C_0\) absolute, such that for
every \(N\ge N_*\) there is a weight \(w_N:[N]\to[1/4,W]\) of
mean one satisfying
\begin{equation}\label{eq:weight-cert}
 I_N(A)\le\delta
 \quad\Longrightarrow\quad
 \frac1N\sum_{x\in A}w_N(x)\le\frac13+\frac{\epsilon}{4}
 \qquad(A\subseteq[N]).
\end{equation}
\end{lemma}

\begin{proof}
Suppose weights are available for all lengths at least \(N_*\), with
height \(W\), threshold \(\delta\) and mass bound \(1/3+e\). At the next
step choose \(\eta\le\eta_0\) and abbreviate \(Q=Q_\eta\),
\(\kappa=\kappa_\eta\), \(B=2Q/\kappa\). For \(N\ge\lceil BN_*\rceil\),
put \(M=\lfloor\kappa N/Q\rfloor\). Then \(M\ge N_*\) and
\(N/B\le M\le\kappa N/Q\). Give every integer weight \(1/4\), and
at \(Qj\) with \(j\in[M]\) add \(3Nw_M(j)/(4M)\). The new weight
has mean one and height at most \(BW\).

For \(A\subseteq[N]\) let \(A_Q=\{j\in[M]:Qj\in A\}\).
Its new weighted mass is
\[
 \frac{|A|}{4N}+\frac{3}{4M}\sum_{j\in A_Q} w_M(j).
\]
Choose the new threshold
\begin{equation}\label{eq:weight-step}
 \delta'=\min\{\delta/B^2,\ \eta\kappa/(144WQ),\ a\eta^2\}.
\end{equation}
Assume \(I_N(A)\le\delta'\). If \(|A|\ge(1/3+\eta)N\),
the third entry in \eqref{eq:weight-step} allows us to use
\eqref{eq:origin}. The mass on the progression is bounded by
\[
 \frac{3}{4M}\sum_{j\in A_Q}w_M(j)
 \le \frac{6W T(A)}{M\eta N}
 \le \frac{12W Q}{\eta\kappa}I_N(A)
 \le \frac1{12},
\]
using \(M\ge\kappa N/(2Q)\) and the second entry in
\eqref{eq:weight-step}. Adding the uniform part gives at most \(1/3\).

If \(|A|<(1/3+\eta)N\), use the old weight bound on the
dilated set \(A_Q\). Dilation gives \(T(A_Q)\le T(A)\),
so \(I_M(A_Q)\le B^2I_N(A)\le\delta\) by the first entry in
\eqref{eq:weight-step}. The new mass is at most
\[
 \tfrac14(1/3+\eta)+\tfrac34(1/3+e)=1/3+e',
 \qquad e'=\tfrac34e+\tfrac14\eta.
\]

Start with weights identically one, \(N_*=W=\delta=1\), \(e_0=2/3\),
and use \(\eta_j=\eta_0(7/8)^j\) at step \(j\).
Inductively \(e_j\le\eta_j/\eta_0\), since
\(3/4+\eta_0/4\le7/8\). Stop at the first \(J\) with
\(\eta_J\le\eta_0\epsilon/4\). This makes the final mass bound
\(1/3+\epsilon/4\), and \(\eta_J>7\eta_0\epsilon/32\).

For the parameter cost, write \(B_j=2Q_{\eta_j}/\kappa_{\eta_j}\)
and \(s_j=\sum_{i<j}\eta_i^{-1}\). Proposition~\ref{prop:shifts}
gives \(\log B_j=O(\eta_j^{-1})\), and the height and minimum length
are multiplied at each step by \(B_j\) and at most \(2B_j\),
respectively. Their logarithms are \(O(s_j)\) at stage \(j\).
By \eqref{eq:weight-step}, the next \(\log(1/\delta')\) is the
maximum of
\[
 \log(1/\delta)+2\log B_j,\qquad
 \log W+\log B_j+\log(72/\eta_j),\qquad
 \log(1/(a\eta_j^2)).
\]
Induction bounds it by a constant times \(s_{j+1}\). Since
\(s_J=O(\eta_J^{-1})=O(\epsilon^{-1})\) and the starting
\(\eta_0\) is absolute, this proves the lemma.
\end{proof}

\subsection{Rounding the weights}

\begin{lemma}\label{lem:round}
Let \(0<\epsilon,\zeta\le1/4\), \(W\ge1\), and \(0\le\theta\le1\).
Suppose \(w:[N]\to[0,W]\) has mean one and
\[
 I_N(S)\le\zeta \quad\Longrightarrow\quad
 \sum_{x\in S}w(x)\le\theta N.
\]
There is an absolute \(C_1\) such that if
\(N\ge C_1\max(2,W,\epsilon^{-1},\zeta^{-1})^{20}\), some
\(B\subseteq[N]\) satisfies
\[
 (1-\epsilon/4)N/W\le |B|\le(1+\epsilon/4)N/W,\qquad
 \alpha(B)\le(\theta+\epsilon)|B|.
\]
\end{lemma}

\begin{proof}
Use the \(3\)-uniform hypergraph \(\mathcal H\) with edges
\(\{x,y,x+y\}\) for \(1\le x<y\) and \(x+y\le N\).
It has \(\lfloor(N-1)^2/4\rfloor\) edges, average degree
\(d\gg N\), maximal pair degree at most three, and maximal
triple degree one. In the notation of Saxton and Thomason
\cite[Definition~3.2]{st}, the co-degree function is therefore
\[
 \delta(\mathcal H,\tau)=4\delta_2+2\delta_3
 \le \frac{12}{d\tau}+\frac{2}{d\tau^2}
 \ll (N\tau^2)^{-1}\quad (0<\tau\le1).
\]
Choose \(\tau=c_1/\sqrt{\zeta N}\) with \(c_1\) a sufficiently
large absolute constant. Then \(\delta(\mathcal H,\tau)\le\zeta/72\)
and \(\tau<1/2\) under the assumed size bound. By
\cite[Corollary~3.6]{st} there is a family \(\mathcal C\) covering
all independent sets with
\[
 e(\mathcal H[S])\le\zeta e(\mathcal H)\quad(S\in\mathcal C),\qquad
 \log|\mathcal C|\ll \sqrt{N/\zeta}\log(1/\zeta)\log N.
\]
A container has at most \(2\zeta e(\mathcal H)+N/2\le\zeta N^2\)
ordered Schur triples for \(N\ge1/\zeta\). The term \(N/2\)
accounts for \(x=y\).

Include each \(x\) in \(B\) independently with probability \(w(x)/W\).
Write \(\mu=N/W\). Each \(S\in\mathcal C\) has
\(\E|B\cap S|\le\theta\mu\). Bernstein's inequality for sums of
independent Bernoulli variables bounds the upper-tail probability of
\(|B\cap S|\) exceeding \((\theta+\epsilon/4)\mu\) by
\(\exp(-c_2\epsilon^2\mu)\) for an absolute \(c_2>0\).
The same bound applies to each tail of \(|B|\) beyond
\(\epsilon\mu/4\) from its mean, since the variance of each sum is
at most \(\mu\).
Put \(R=\max(2,W,\epsilon^{-1},\zeta^{-1})\). We have
\(\epsilon^2\mu\ge N/R^3\) and
\(\log|\mathcal C|\ll N^{3/4}R^{3/2}\), using
\(\log N\ll N^{1/4}\) and \(\log(1/\zeta)\le R\).
The ratio \((N/R^3)/(N^{3/4}R^{3/2})=N^{1/4}/R^{9/2}\) is at
least \(C_1^{1/4}R^{1/2}\) for \(N\ge C_1R^{20}\). Taking \(C_1\)
large makes \((|\mathcal C|+2)\exp(-c_2\epsilon^2\mu)<1\). In a successful sample every sum-free subset of \(B\)
is contained in a container, so its relative size is at most
\((\theta+\epsilon/4)/(1-\epsilon/4)\le\theta+\epsilon\).
\end{proof}

\begin{proof}[Proof of Theorem~\ref{thm:main} assuming Proposition~\ref{prop:shifts}]
Apply Lemma~\ref{lem:weights} with \(0<\epsilon\le1/12\) and
Lemma~\ref{lem:round} with
\(\theta=1/3+\epsilon/4\), \(\zeta=\min(\delta,1/4)\).
For every \(N\ge\exp(C_2/\epsilon)\) the rounding succeeds for an
absolute \(C_2\ge C_0\), by the polynomial size condition.
For an integer \(n\ge\exp(C_2/\epsilon)\), take
\(N=\lceil(1+\epsilon)Wn\rceil\). Enlarging \(C_2\) so that
\(n\ge1/\epsilon\), we have
\((1+\epsilon)Wn\le N\le(1+2\epsilon)Wn\).
The size of the rounded set satisfies \(n<|B|\le(1+3\epsilon)n\).
Choose an arbitrary \(n\)-element subset \(A\) of \(B\). Then
\[
 \alpha(A)\le(1/3+5\epsilon/4)(1+3\epsilon)n
             \le n/3+3\epsilon n,\qquad
 \max A\le 2n\exp(C_0/\epsilon).
\]
For large \(n\) put \(\epsilon=2C_2/\log n\). Its size threshold
is \(\sqrt n\le n\), and the bound on the largest element is
at most \(2n^{3/2}\le n^2\). This proves the theorem.
\end{proof}

\section{Sampling a cyclic profile}\label{sec:sampling}

For a finite abelian group \(G\), write \(m\) for its uniform
probability measure. For a measurable \(h:G\times(0,1]\to[0,1]\),
put \(h_s(g)=h(g,s)\) and \(\Delta=\{(u,w): u,w>0,\ u+w\le1\}\).
Define
\[
 d(h)=\int_0^1\E_g h_s(g)\,ds,\qquad
 J(h)=\int_\Delta \E_{g,g'}h_u(g)h_w(g')h_{u+w}(g+g')\,du\,dw .
\]
If \(h\) is an indicator, we identify its section \(h_s\) with a
subset of \(G\) when writing a measure or symmetric difference.
A set is independent under a shift \(v\) if it is disjoint from its
translate by \(v\).

Sampling introduces a shifted sum in the integer count. For
\(F\subseteq[N]\) and \(j\in\Z\), write
\[
 T_j(F)=\sum_{x,y\in\Z}\ind_F(x)\ind_F(y)\ind_F(x+y+j)
\]
for the count with output \(x+y+j\), so \(T_0(F)=T(F)\).

\begin{lemma}[A shift with a small translated triple count]\label{lem:return}
Let \(F\subseteq[N]\) with \(|F|\ge N/3\), let \(b,t\) be positive
integers, and let \(0<\xi<1\). Put
\(k=\lceil32/\xi^2\rceil\) and \(L=4\cdot6^k\).
If \(btL\le N\), there is a multiple \(q\) of \(bt\) with \(bt\le q<btL\)
such that
\[
 T_{\pm q}(F)\le 2T(F)+\xi^2N^2.
\]
\end{lemma}

\begin{proof}
We use a sampling and translation argument of Croot and Sisask
\cite{croot-sisask}, with the square root of the convolution as the
quantity to approximate. Put \(f=\ind_F\), \(U=f*f\) in counting
measure on \(\Z\), \(a=|F|\), and \(P=U/a\). For a tuple
\(\mathbf u=(u_1,\dots,u_k)\) sampled independently and uniformly
from \(F\), put \(P_{\mathbf u}(x)=k^{-1}\sum_i f(x-u_i)\).
At \(P(x)>0\),
\[
 \E\big(\sqrt{P_{\mathbf u}(x)}-\sqrt{P(x)}\big)^2
 \le \frac{\E(P_{\mathbf u}(x)-P(x))^2}{P(x)}
 =\frac{1-P(x)}{k}.
\]
Both terms are zero where \(P=0\). Thus at least half the tuples have
\(\|\sqrt{P_{\mathbf u}}-\sqrt P\|_2^2\le 2V/k\), with
\(V=|F+F|\le2N\).

Translate each of these tuples by \(s\) in all coordinates for
each \(s\in S=bt\{0,\dots,L-1\}\). There are at least \(a^kL/2\)
tuple-shift pairs and at most \(|F+S|^k\le(6a)^k\) distinct
images, since \(btL\le N\). Some image \(\mathbf z\) has two
preimages with distinct shifts \(s,s'\). For each such shift,
\(P_{\mathbf z}(\cdot+s)=P_{\mathbf z-s\mathbf 1}\).
Comparing the two approximations after translating one of them,
and taking \(q=|s'-s|\), gives
\[
 \|\sqrt{U(\cdot\pm q)}-\sqrt U\|_2^2\le 8aV/k.
\]
The norms for the two signs are equal by translation. Since
\(T_j(F)=\|\sqrt{U(\cdot-j)}\sqrt f\|_2^2\), the triangle
inequality and \(0\le f\le1\) give
\[
 T_{\pm q}(F)\le 2T(F)+16aV/k
              \le 2T(F)+32N^2/k.
\]
The choice of \(k\) finishes the proof.
\end{proof}

\begin{samepage}
\begin{lemma}[Sampling modulo \(q\)]\label{lem:sampler}
Let \(F\subseteq[N]\) be independent under the shift \(t\), where
\(1\le t\le q\le N\). In \(G=\Z/q\Z\) define
\[
 n_r(s)=r+q\left\lfloor\frac{Ns-r}{q}\right\rfloor,\qquad
 h(r,s)=\ind_F(n_r(s))\,\ind_{\{Ns-n_r(s)\ge t\}}.
\]
Then each section is independent under \(v=t\bmod q\), and
\begin{equation}\label{eq:sampler}
 d(h)\ge |F|/N-q/N-t/q,\qquad J(h)\le (T(F)+T_q(F))/N^2.
\end{equation}
\end{lemma}
\end{samepage}

\begin{proof}
The integer \(n_r(s)\) is the largest one in residue \(r\) at most
\(Ns\), independent of the representative of \(r\). If \(h(r,s)=1\),
then \(n_r(s)+t\) is the largest such integer in residue \(r+t\)
and is outside \(F\). For fixed \(s\), the deficits
\(Ns-n_r(s)\) run through \(\{Ns\}+\{0,\ldots,q-1\}\), where
\(\{Ns\}\) is the fractional part. The cutoff deletes exactly \(t\)
residues. Before deletion, an integer \(n\in F\) is sampled only
in its residue, on \([n/N,(n+q)/N)\cap(0,1]\). Its cell has
length \(q/N\) when \(n\le N-q\). This proves the density bound.

At a time pair \((u,w)\in\Delta\), write \(a=n_r(u)\),
\(b=n_{r'}(w)\). Each deficit is in \([0,q)\), so their sum
is in \([0,2q)\). The largest integer at most \(N(u+w)\) in
residue \(r+r'\) is therefore \(a+b\) or \(a+b+q\).
Thus the integrand is at most
\(\ind_F(a)\ind_F(b)(\ind_F(a+b)+\ind_F(a+b+q))\).
The two inputs fix the residue pair and their cell rectangle
has area at most \((q/N)^2\). After the uniform residue average
their contribution is at most \(1/N^2\) times the displayed
bound. Summing over integer pairs proves the count estimate.
\end{proof}

\section{Structure of a cyclic profile at fixed accuracy}\label{sec:structure}

For integer \(D\ge1\) let \(Q(D)=\lcm(1,\dots,D)\). A character is
an additive homomorphism to \(\T\). A character of order \(\ell\) has
image \(C_\ell=(1/\ell)\Z/\Z\). An interval in \(C_\ell\) consists of
consecutive residues in cyclic order.

\begin{proposition}\label{prop:structure}
For every \(0<\varepsilon\le1/20\) there are an integer \(D\ge1\) and
a number \(J_*>0\) with the following property. Suppose \(G\) is
finite cyclic, \(v\in Q(D)G\), and \(h\) is an indicator on
\(G\times(0,1]\) with \(v\)-independent sections almost everywhere,
\(d(h)\ge1/3\), and \(J(h)\le J_*\).
There are a character \(\chi:G\to\T\) of order \(\ell>D\) and a
real \(\lambda\) such that
\begin{equation}\label{eq:profile}
 \int_0^1\E_g\big|h(g,s)-\ind_{[1/3,2/3)}(\chi(g)-\lambda s)\big|\,ds
 \le\varepsilon,\qquad \|\chi(v)\|_{\T}\ge 5/16.
\end{equation}
Here \(\|\cdot\|_{\T}\) is distance to the nearest integer.
\end{proposition}

In the application, the shift is \(v=t\bmod q\).
Divisibility of \(t\) by \(Q(D)\) ensures the shift hypothesis,
and the resulting character has order \(\ell\mid q\).
The structure argument proceeds as follows.
Most slices have density near \(1/3\) (\S\ref{sec:slice-densities}),
so a typical pair has interval models under a common character
(\S\ref{sec:critical-pairs}). Fixing one time then gives a single
character for almost every slice, by uniqueness up to sign
(\S\ref{sec:common-character}). The remaining freedom is the
position of each interval. The circle triple count bounds the
failure of the positions to add, yielding a linear phase in time
(\S\ref{sec:linear-phase}).

\subsection{Section densities}\label{sec:slice-densities}

\begin{lemma}\label{lem:finite-exclusion}
Let \(X,Y,Z\) be subsets of a finite abelian group \(G\), each
independent under \(v\). Suppose their total normalized size is
at least \(1+\sigma\), with \(0<\sigma\le1/2\), and their normalized
cross count \(\E_{x,y}\ind_X(x)\ind_Y(y)\ind_Z(x+y)\) is at most
\(\sigma^2/4096\). Some subgroup of index at most \(2/\sigma\)
does not contain \(v\).
\end{lemma}

\begin{proof}
Write \(M=|G|\), \(U=\sigma M\), and \(E\) for the unnormalized
cross count. Each set has size at most \(M/2\), hence at least \(U\).
If \(U<32\) then \(E=0\). Otherwise take \(j=\lfloor U/16\rfloor\),
so \(U/32\le j\le U/16\) and \(1\le j\le\min(|Y|,|Z|)\).
Let \(r(x)\) be the representation count in \(Z-Y\). Then
\[
 \sum_x\min(r(x),j)\le j(M-|X|)+E
 \le j(|Y|+|Z|)-jU+E .
\]
The last bound is strictly below \(j(|Y|+|Z|)-2j^2+1\),
since \(jU-2j^2-E\ge (1/32-1/128-1/4096)U^2>0\).
This is the deletion case of Grynkiewicz's popular-sum theorem
\cite[Theorem~1.1]{grynkiewicz}. It gives subsets of \(Z,Y\)
with at most \(j-1\) total deletions whose difference set is
exactly \(\{x:r(x)\ge j\}\).
Remove its intersection with \(X\), of size at most \(E/j\le U/128\).

In either case the retained sets \(X',Y',Z'\) have no cross
triple, each is nonempty, and their total size \(M'\) is at least
\(M+U/2\). Let \(H\) be the stabilizer of \(Z'-Y'\).
The saturation \(X'+H\) is disjoint from \(Z'-Y'\).
Kneser's theorem \cite[Theorem~1]{devos} gives
\[
 |X'+H|+|Y'+H|+|Z'+H|\le M+|H|.
\]
Thus \(|H|\ge U/2\). If \(v\in H\), each saturation contains
both its retained set and its disjoint translate. This would give
\(2M'\le M+|H|\le2M\), a contradiction.
\end{proof}

To relate the total density in the three slices to the average
density of the profile, use the probability measure
\(d\mu=3(u+w)\,du\,dw\) on \(\Delta\).
At input time \(s\), integrating \(3(s+w)\) over \(0<w<1-s\)
gives density \(3(1-s^2)/2\). At sum time \(s\), the change of
variables \((u,w)\mapsto(u,u+w)\) gives density \(3s^2\).
The three marginal densities add to three. The weight itself
is at most three on \(\Delta\).

\begin{lemma}\label{lem:triangle}
For \(g\in L^1(0,1)\) and \(R=\int_\Delta|g(u)+g(w)+g(u+w)|\,d\mu\),
one has \(\|g\|_1\le4500 R\).
\end{lemma}

A proof is given in Appendix~\ref{app:slice-density}.

Let \(d_s=m(h_s)\). Each section has \(d_s\le1/2\). For a time pair,
write \(s=u+w\), \(S=d_u+d_w+d_s-1\), and \(e(u,w)\) for the normalized
spatial count in \(J(h)\). Since \(S\le1/2\), if \(S\ge2/D\)
and \(e(u,w)\le S^2/4096\),
Lemma~\ref{lem:finite-exclusion} gives a subgroup avoiding \(v\)
whose index is at most \(D\). This conflicts with \(v\in Q(D)G\),
since the index divides \(Q(D)\). Hence \(S_+^2\le (2/D)^2+4096e(u,w)\)
almost everywhere. The mean of \(S\) under \(\mu\) is
\(3d(h)-1\ge0\), so \(\E_\mu|S|\le2\E_\mu S_+\). Cauchy--Schwarz
and Lemma~\ref{lem:triangle} applied to \(d_s-1/3\) give
\begin{equation}\label{eq:flat}
 M:=\int_0^1 |d_s-1/3|\,ds
 \le9000\sqrt{(2/D)^2+12288J(h)}.
\end{equation}

\subsection{Near-critical cyclic pairs}\label{sec:critical-pairs}

\begin{lemma}\label{lem:pairs}
Given \(0<e_0\le1/100\) and an integer \(D_1\), there exist
\(D_0\ge\max(1,D_1)\) and \(\nu,E_0>0\) as follows.
Let \(D\ge D_0\), let \(G\) be finite cyclic and \(v\in Q(D)G\).
Suppose \(A,B,C\subseteq G\) each have mass within \(\nu\) of
\(1/3\), the inputs \(A,B\) are independent under \(v\), and their
normalized cross count with \(C\) is at most \(E_0\).
Then \(A,B\) admit interval preimages \(P_A,P_B\) under one
character of order exceeding \(D\), with
\(m(A\triangle P_A)+m(B\triangle P_B)<e_0\).
\end{lemma}

\begin{proof}
Use Griesmer's near-critical sumset theorem
\cite[Theorem~1.15]{griesmer} at precision \(e_2=e_0/10\).
It gives a subgroup index bound \(D_2\) and a positive sumset
tolerance. Take \(D_0=\max(D_2,D_1,1)\) and set
\(e_1=\min(e_2,1/(100D_0))\). Let \(d_1>0\) be the
tolerance in the popular-sumset theorem \cite[Theorem~1.1]{griesmer}
at precision \(e_1\). Choose \(\nu=\min(e_0/10,d_1/6)\) and
\(E_0=d_1^2/2\).

Let \(a,b,c\) be the masses of \(A,B,C\).
The set \(\{z:(\ind_A*\ind_B)(z)>d_1\}\), with normalized
convolution, has mass at most \(1-c+E_0/d_1\le a+b+d_1\).
Both input masses exceed \(e_1\). Thus the popular-sumset theorem
gives sets \(S,T\) with
\(m(S+T)\le s+t\), where \(s=m(S)\), \(t=m(T)\), and
\(e_A+e_B<e_1\) for \(e_A=m(A\triangle S)\),
\(e_B=m(B\triangle T)\). Both \(s,t>e_2\). Apply the near-critical
sumset theorem to \(S,T\), obtaining a subgroup \(K\) of index
at most \(D_2\). It contains \(v\), and \(m(K)\ge1/D_0\).

The overlaps of \(S\) and \(T\) with their translates by \(v\)
have masses at most \(2e_A,2e_B\), respectively, by independence
of \(A,B\). Put \(U=S\cap A\), \(V=T\cap B\).
These sets are independent under \(v\).

\paragraph{Quasi-periodic case.}
Griesmer's case (II) gives a nonempty subset \(S_1\subseteq S\)
or \(T_1\subseteq T\) whose \(K\)-saturation adds at most
\(e_2m(K)\) to its mass. If the subset is \(S_1\), then \(S_1\)
and \(S_1+v\) lie in the same saturation and overlap in mass
at least \((1-2e_2)m(K)>2e_1\). The overlap bound rules this out,
and the same argument rules out \(T_1\).

\paragraph{Periodic sum case.}
In case (I), \(P=S+T+K\) has mass at most
\(s+t+e_2m(K)<1\). Let \(H\) be its stabilizer, \(h=m(H)\le1/2\),
and \(c'=1-m(P)\). Set \(p_S=m(S+H)\), \(p_T=m(T+H)\).
By Kneser's theorem applied to \((S+K)+(T+K)\),
\(p_S+p_T+c'\le1+h\). The shift belongs to \(H\), so
\(2m(U)\le p_S\) and \(2m(V)\le p_T\). Also \(p_S\ge h\).
It follows that
\[
 m(U)+m(V)+c'
 \le 1+h-\tfrac12(p_S+p_T)\le 1+h/2-m(V)
 \le 11/12+\nu+e_1 .
\]
On the other hand \(m(U)+m(V)+c'\ge s+t-e_1+1-m(P)
\ge1-e_1-e_2\), a contradiction since
\(\nu+2e_1+e_2<1/12\).

\paragraph{Interval case.}
In case (III), both sets lie in translates of parallel intervals
or progressions in \(K\) that add less than \(e_2\) to each mass.
In particular \(T\) lies in one coset of \(K\). Independence of
\(V\) under \(v\in K\) gives \(m(K)\ge2m(V)>1/2\), hence \(K=G\).
A finite group has no circle quotient. Thus the parallel
progressions supplied by the theorem are full preimages of intervals
under one cyclic quotient character. Absorb the translates into the
intervals, giving the desired \(P_A,P_B\) with error less than
\(e_1+2e_2\).

If the order \(\ell\) were at most \(D\), then \(\ell\mid Q(D)\) and the
character would kill \(v\). An invariant \(P_A\) can contain at most
half its mass from the independent set \(A\), giving
\(m(A\triangle P_A)\ge a>e_0\). This proves the order bound.
\end{proof}

\subsection{One character for the slices}\label{sec:common-character}

We will use the following uniqueness observation at a fixed density
margin. Put \(b_{\mathrm F}=\sin(\pi/4)/(2\pi)\) and
\(D_1=\lceil b_{\mathrm F}^{-2}+1\rceil\). Suppose \(P,Q\) are interval
preimages under characters \(\chi,\psi\) of orders above \(D_1\),
both lifts have densities in \([1/4,3/4]\), and
\(m(P\triangle Q)<b_{\mathrm F}\). Then \(\chi=\pm\psi\).

Indeed, the Fourier coefficients
\(\widehat{\ind_P}(\gamma)=\E_g\ind_P(g)e^{-2\pi i\gamma(g)}\)
vanish outside multiples of \(\chi\) by invariance under
\(\ker\chi\). For a lift of \(k\) consecutive residues under
a character \(\theta\) of order \(\ell\), the geometric series gives
coefficient magnitude
\[
 \frac{|\sin(\pi j k/\ell)|}{\ell|\sin(\pi j/\ell)|}
 \quad (0<|j|\le \ell/2)
\]
at \(j\theta\). At \(j=1\) the density margin and
\(\ell\sin(\pi/\ell)\le\pi\) give a lower bound of \(2b_{\mathrm F}\).
For nonzero signed \(j\), the chord bound
\(\sin(\pi |j|/\ell)\ge2|j|/\ell\) gives the upper bound \(1/(2|j|)\).
The coefficients of the two lifts differ by less than \(b_{\mathrm F}\).
Thus \(\chi=j\psi\) and \(\psi=i\chi\) for signed integers
\(0<|i|,|j|<1/(2b_{\mathrm F})\). The relations give the same order \(\ell\)
for both characters and \(\ell\mid ij-1\). Since \(|ij-1|<D_1<\ell\),
we have \(ij=1\), proving uniqueness up to reflection.

\subsection{Recovering a linear phase}\label{sec:linear-phase}

\begin{lemma}\label{lem:phase}
For every \(\epsilon>0\) there is \(\delta_{\mathrm p}>0\)
such that any measurable \(f:(0,1]\to\mathbb C\) with \(|f|=1\)
almost everywhere and
\[
 R=\int_\Delta |f(u)f(w)-f(u+w)|\,du\,dw\le\delta_{\mathrm p}
\]
satisfies \(\int_0^1|f(s)-e^{2\pi i\lambda s}|\,ds<\epsilon\)
for some \(\lambda\in\R\).
\end{lemma}

\begin{proof}
It suffices to take \(0<\epsilon\le1\).
Put \(J_0=(0,1/2)\), \(a=\epsilon/16\), \(h_0=\epsilon^2/10000\).
Apply Christ and Shao's local approximate-character theorem
\cite[Proposition~5.1]{christ-shao} on \(J_0\) at \(L^2\) accuracy
\(h_0\). It has a tolerance \(d_*>0\) for which the product of two
phase values must differ from a function of their sum by more than
\(d_*\) on a set of area less than \(d_*\).
Use a measurable argument of \(f\) for the phase and \(f\) for
the sum function. With \(\delta_{\mathrm p}=\min(h_0,d_*^2/2)\),
Markov's inequality gives exceptional area at most \(R/d_*\le d_*/2\)
on \(J_0^2\). The conclusion is a function \(ce^{i\xi s}\) with
\(|c|=1\) and \(L^1(J_0)\) error \(k\le h_0\), by Cauchy--Schwarz.

The sum of the two inputs on \(J_0^2\) has unnormalized density
\(\ell(s)=\min(s,1-s)\). Replacing the two inputs by \(ce^{i\xi s}\)
and using the residual gives
\[
 \int_0^1\ell(s)|f(s)-c^2 e^{i\xi s}|\,ds\le R+k .
\]
On \([1/4,1/2]\), compare the two approximations using \(\ell\ge1/4\).
This gives \(|c-c^2|/4\le k+4(R+k)\) and
\(|1-c^2|\le32R+40k\).
On \([a,1-a]\) we have \(\ell\ge a\), and outside that interval
the unweighted integral costs at most \(4a\).
Thus the \(L^1(0,1)\) error from \(e^{i\xi s}\) is bounded by
\(4a+(R+k)/a+32R+40k<\epsilon\). Take \(\lambda=\xi/(2\pi)\).
\end{proof}

We record how a circle triple count measures phase error. For
arcs \(I_s=[1/3,2/3)+p(s)\) and a time pair \((u,w)\), put
\(d=\|p(u)+p(w)-p(u+w)\|_{\T}\).
Write \(E_{\T}(u,w)\) for the normalized circle triple count of
\(I_u,I_w,I_{u+w}\) in their indicated roles. The support of the
input convolution is an arc of length \(2/3\). Its complement has
the same length as the output arc, whose center differs from the
complement's center by distance \(d\). Translating the support to
\([0,2/3]\) gives convolution density
\(C(x)=\min(x,2/3-x)\) there, by taking the length of
\([0,1/3]\cap[x-1/3,x]\). Reflection permits taking the signed
shift to be \(d\ge0\). With \(H(t)=\int_0^t C(x)\,dx\),
the count on the complement of the shifted support is
\[
 E_{\T}(u,w)=H(d)-H((d-1/3)_+).
\]
The second term removes the part where the shifted support wraps
around zero. For \(d\le1/3\) this is \(d^2/2\).
For \(d>1/3\) its derivative in \(d\) is \(1-2d\ge0\).
Consequently
\begin{equation}\label{eq:circle-count}
 E_{\T}(u,w)\ge \tfrac12\min(d,1/3)^2\ge \tfrac29 d^2 .
\end{equation}

\begin{proof}[Proof of Proposition~\ref{prop:structure}]
Given \(\varepsilon\), take the tolerance \(\delta_{\mathrm p}\)
in Lemma~\ref{lem:phase} at error \(\varepsilon/2\).
Choose \(0<\rho\le1/100\) such that
\[
 3\pi\sqrt{15\rho}\le\delta_{\mathrm p},\qquad
 9\rho+\varepsilon/4\le\varepsilon,\qquad 12\rho\le1/48 .
\]
Put \(\tau=\rho/4\),
\(e_0=\min(\rho/4,b_{\mathrm F}/4,1/100)\), and take \(D_0,\nu,E_0\)
from Lemma~\ref{lem:pairs} with \(D_1\) as above.
Choose \(D\ge\max(D_0,\rho^{-1})\) and \(0<J_*\le\rho\) so that
\begin{equation}\label{eq:bad-times}
 M_*:=9000\sqrt{(2/D)^2+12288J_*}\le\rho,\qquad
 3M_*/\nu+3J_*/E_0\le\tau^3/100 .
\end{equation}
Taking \(D\) large and \(J_*\) small makes both inequalities hold,
with \(\nu,E_0,\tau\) fixed. By \eqref{eq:flat}, \(M\le M_*\).

The left-hand side of \eqref{eq:bad-times} bounds the
\(\mu\)-measure of the pairs for which one of the three masses
is farther than \(\nu\) from \(1/3\), or the spatial count
exceeds \(E_0\). Indeed, the marginal densities add to three
and the count is nonnegative. The time marginals are absolutely
continuous, so one may also discard times where shift independence
fails. On each remaining pair Lemma~\ref{lem:pairs} gives parallel
interval lifts of the two inputs, with total error less than \(e_0\)
and order greater than \(D\). Their densities lie in \([1/4,3/4]\)
by the bounds on \(\nu,e_0\). The choices can be taken measurably
by ordering the finitely many possibilities.

\paragraph{Fixing a reference slice.}
For \(w\in[\tau/2,\tau]\), the rectangle of times
\(u\in(0,1-\tau]\) lies in \(\Delta\) and \(\mu\) has density
at least \(3\tau/2\) there. The average length of exceptional
\(u\)'s over this interval is at most
\(4(\tau^3/100)/(3\tau^2)<\tau\). Choose \(w_0\)
with exceptional length less than \(\tau\) and a measurable slice.
This leaves a good \(u_0\). Compare the lift of \(h_{w_0}\)
from \((u_0,w_0)\) with its lift from each good \((u,w_0)\).
They differ by at most \(2e_0<b_{\mathrm F}\), so the characters agree up
to sign. Reflect the interval for \(h_u\) if necessary to use
the character \(\chi\) from \((u_0,w_0)\).
Take the empty interval on exceptional times. The measurable
interval lifts \(P_s\) under \(\chi\) then satisfy
\[
 \delta:=\int_0^1 m(h_s\triangle P_s)\,ds\le e_0+2\tau\le\rho .
\]

\paragraph{Passing to the circle.}
Write \(\ell>D\) for the order of \(\chi\). Represent the quotient
interval of each \(P_s\) by thickening its residues with cells of
width \(1/\ell\) centered at them. Choose a center for the empty
interval as well. Change the width to \(1/3\), keeping the center
and using half-open endpoints to obtain a circle arc \(I_s\).
The two concentric arcs are nested, and their grid counts differ
from their lengths times \(\ell\) by at most one each. Therefore,
for \(B_s=\chi^{-1}(I_s)\),
\[
 \int_0^1 m(h_s\triangle B_s)\,ds\le
 a:=2\delta+M+2/\ell\le 5\rho .
\]
In this estimate the discrepancy between the length of the
first arc and \(m(h_s)\) is at most the original lift error at
time \(s\).

Replacing \(h_s\) by \(B_s\) in the triple count costs at most \(2a\)
because the two unweighted input marginals have density \(1-s\)
and the output marginal has density \(s\). The spatial average for
\(B_s\) uses uniform points of \(C_\ell\). Add independent uniform
displacements in \([0,1/\ell)\) to the two grid inputs to get
independent uniform circle points. Each input test changes with
probability at most \(2/\ell\), from the two arc endpoints. The
output grid point is uniform conditional on the displacements
and moves by less than \(2/\ell\), so its test changes with
probability at most \(4/\ell\). The triangle has area \(1/2\).
Thus the integrated circle count of the \(I_s\) is at most
\[
 \int_\Delta E_{\T}(u,w)\,du\,dw\le J(h)+2a+4/\ell\le15\rho .
\]
\paragraph{Aligning the phases.}
Write \(I_s=[1/3,2/3)+p(s)\). By \eqref{eq:circle-count} and
Cauchy--Schwarz on the triangle, the residual of \(e^{2\pi i p(s)}\)
is at most \(3\pi\sqrt{15\rho}\), since the chord at distance \(d\)
has length at most \(2\pi d\). Lemma~\ref{lem:phase} thus gives
an approximation \(e^{2\pi i\lambda s}\) of \(L^1\) error
less than \(\varepsilon/2\). Moving an arc a circular distance
\(d\) changes at most \(2d\) of its mass, at most half the chord.
Its grid symmetric difference adds at most \(4/\ell\) to this bound,
from the endpoints of the two arcs. Thus the approximation in
\eqref{eq:profile} has error at most \(a+\varepsilon/4+4/\ell
\le9\rho+\varepsilon/4\).

Finally, the shift independence of \(h_s\) gives
\(m(B_s\cap(B_s+v))\le 2m(h_s\triangle B_s)\).
The circle overlap of \(I_s\) with its translate by \(\chi(v)\)
has length \((1/3-\|\chi(v)\|_\T)_+\) and a grid discrepancy
at most \(2/\ell\). Integrating gives
\(\|\chi(v)\|_\T\ge1/3-2a-2/\ell\ge5/16\).
\end{proof}

\section{Returning to the integers}\label{sec:return-integers}

Write \(I=[1/3,2/3)\) for the middle third in \(\T\).
The slope in the profile becomes part of an integer rotation
frequency. Controlling the phase displacement at the prescribed
shift will give discrepancy estimates on every block.

\begin{lemma}\label{lem:unfold}
Suppose \(h\) is the profile in Lemma~\ref{lem:sampler},
\(t/q,q/N\le1/10\), and its error in \eqref{eq:profile} is at most
\(\varepsilon\le1/20\), with \(\chi(r)=kr/q\bmod1\).
Put \(\omega=k/q-\lambda/N\). Some \(\beta\in\R\) gives
\begin{equation}\label{eq:affine}
 |F\triangle\{n\in[N]:\omega n+\beta\bmod1\in I\}|/N
 \le \varepsilon/(1-t/q)+q/N .
\end{equation}
Moreover
\begin{equation}\label{eq:unfold-displacement}
 |\lambda|(q-t)/N<4,\qquad
 \|\omega t-\chi(t\bmod q)\|_{\T}<4t/(q-t).
\end{equation}
\end{lemma}

\begin{proof}
For \(n\le N-q\) in residue \(r\), the retained sampling cell
at time \(s=n/N+\sigma\) has \(t/N\le\sigma<q/N\).
On it the sampled value is \(\ind_F(n)\) and the model is
\(\ind_I(\omega n-\lambda\sigma)\).
The cells are disjoint within each residue. Restricting the
approximation integral to them gives
\[
 \frac1q\int_{t/N}^{q/N}
 \sum_{n=1}^{N-q}|\ind_F(n)-\ind_I(\omega n-\lambda\sigma)|
 \,d\sigma\le\varepsilon .
\]
Some \(\sigma\) has at most \(\varepsilon qN/(q-t)\)
mismatches on the first \(N-q\) integers. This gives
\eqref{eq:affine} with \(\beta=-\lambda\sigma\).

If the phase traversed an interval of length
\(L=|\lambda|(q-t)/N\ge4\) on each cell, the proportion in \(I\)
would be at least \(1/3-1/L\), and the proportion outside at
least \(2/3-1/L\), by splitting into unit intervals and a remainder.
Each cell would then disagree with its constant sampled value on
at least a twelfth of its length. The integral would be at least
\((1-q/N)(1-t/q)/12\ge27/400>1/20\).
This proves the slope bound, which yields the displacement bound
since \(\omega t-kt/q=-\lambda t/N\).
\end{proof}

\begin{lemma}\label{lem:discrepancy}
Let \(H,t\) be positive integers with \(Q(H)\mid t\).
If \(b=\|\omega t\|_\T>0\), then for every consecutive block \(J\)
of \(M\) integers, real offset \(\beta\), and circle interval \(V\)
with arbitrary endpoint conventions,
\begin{equation}\label{eq:discrepancy}
 \big|\#\{n\in J:\omega n+\beta\bmod1\in V\}-M|V|\big|
 \le \frac{M}{H+1}+\frac{3t}{2b}.
\end{equation}
\end{lemma}

\begin{proof}
Since \(h\mid t\), we have \(\|h\omega\|_\T\ge hb/t\)
for \(1\le h\le H\).
Use the constant-target case of the Erd\H{o}s--Tur\'an
inequality of Martin and Sitar \cite[Theorem~1.2]{martin-sitar}.
For a closed arc the bound is \(M/(H+1)\) plus
\[
 \sum_{h=1}^H \left(\frac{2-c}{H+1}+\frac{c}{h}\right)
  \max_{m\le M}\left|\sum_{j=1}^m e^{2\pi i h(\omega j+\beta')}\right|,
 \qquad c=16/(7\pi),
\]
where \(\beta'\) accounts for the start of the block.
The geometric series bounds the maximum by
\(1/(2\|h\omega\|_\T)\le t/(2bh)\).
The coefficient multiplying \(t/(2b)\) is at most
\((2-c)H/(H+1)+2c<3\), using \(\sum_{h\ge1}h^{-2}\le2\).
For other endpoints use the closed arc above for the upper
bound and the closed complement above for the lower bound.
\end{proof}

\begin{lemma}[Removing the phase offset]\label{lem:offset}
Let \(N\ge16\), \(0<\delta\le1/100\), and \(\gamma=\delta/100\).
Write \(\phi(n)=\omega n\bmod1\) and
\(B_\beta=\{n\in[N]:\phi(n)+\beta\in I\}\).
Suppose the discrepancy of \(\phi\) on any consecutive block
of \([N]\) and any interval is at most \(\mathcal D\), with
\(\mathcal D/N\le\gamma/64\). If \(F\subseteq[N]\) satisfies
\[
 |F\triangle B_\beta|/N\le\gamma^2/4096,\qquad
 I_N(F)\le\gamma^2/1024,
\]
then there is \(\overline F\subseteq[N]\) such that
\begin{align}
 &\{n:\phi(n)\in[1/3+3\delta,2/3-3\delta]\}
      \ \subseteq\ \overline F\
      \subseteq\ \{n:\phi(n)\in[1/3-\delta,2/3+\delta]\}, \label{eq:strip}\\
 &T(\overline F)\le T(F),\qquad
 |\overline F|\ge |F|-\frac{8T(F)}{\delta N}. \notag
\end{align}
The sets in \eqref{eq:strip} are restricted to \([N]\).
\end{lemma}

\begin{proof}
Let \(C=\ind_I*\ind_I\) on the circle with Haar measure.
Its value decreases linearly from \(1/3\) at zero to zero at
\(1/3\), vanishes through \(2/3\), and increases linearly to
\(1/3\) at one. This follows from the length of \(I\cap(\theta-I)\)
at phase \(\theta\).
Write \(\theta_z=\phi(z)+\beta\) and let \(r_{B_\beta}(z)\)
count ordered representations of \(z\) as a sum from \(B_\beta\).
The input phases for such representations lie in
\(I\cap(\theta_z+\beta-I)\). This region is a union of at most
two arcs of total length \(C(\theta_z+\beta)\).
Using the prefix \(1\le x<z\) gives
\begin{equation}\label{eq:affine-convolution}
 r_{B_\beta}(z)\ge (z-1)C(\theta_z+\beta)-2\mathcal D .
\end{equation}
Changing a set at \(j\) indices changes the triple count by
at most \(3Nj\). In each role a fixed changed index leaves at
most \(N\) choices. Thus \(I_N(B_\beta)\le7\gamma^2/4096\).

Choose \(\beta\in[-1/2,1/2]\). If \(\beta\ge\gamma\),
then \(C(\theta+\beta)\ge\gamma/2\) on the arc
\([2/3-\gamma/2,2/3-\gamma/4]\subset I\).
Indeed, before reduction the shifted arc is in
\([2/3+\gamma/2,7/6-\gamma/4]\), and the claim holds on both
linear branches of \(C\). For \(\beta\le-\gamma\), reflect to
the arc \([1/3+\gamma/4,1/3+\gamma/2]\) using \(C(-x)=C(x)\).
In either case, at least \(N\gamma/8-\mathcal D\ge N\gamma/16\)
indices \(z>N/2\) have \(\theta_z\) in that arc. For each,
\eqref{eq:affine-convolution} is at least \(N\gamma/8\).
This would give \(I_N(B_\beta)\ge\gamma^2/128\).
It follows that \(\|\beta\|_\T<\gamma\), and the symmetric
difference between \(I\) and \(I-\beta\) gives
\[
 |F\triangle B_0|/N\le\gamma^2/4096+2\gamma+2\mathcal D/N
 <\delta/32.
\]

Put \(S=\{n\in[N]:\phi(n)\in[1/3-\delta,2/3+\delta]\}\).
Outside this arc \(C(\phi(n))\ge\delta\).
If \(z\in F\setminus S\) and \(z>N/2\), then
\eqref{eq:affine-convolution} at \(\beta=0\) gives at least
\((z-1)\delta-2\mathcal D\) baseline representations.
Changing both inputs to \(F\) costs at most \(2|F\triangle B_0|\).
At least \(\delta N/4\) representations remain, by \(N\ge16\)
and the discrepancy and symmetric-difference bounds.
If \(x\in F\setminus S\) and \(x\le N/2\), instead count
partners \(y\le N-x\) with \(y,x+y\in B_0\). The admissible
phases lie in \(I\cap(I-\phi(x))\), of length \(C(\phi(x))\)
since \(I=-I\) up to endpoints. The count is at least
\((N-x)\delta-2\mathcal D\) before changing the partner and
output indicators to \(F\), and at least \(\delta N/4\) after.
Summing in the output role for the late indices and in the
input role for the early ones gives \(|F\setminus S|\le8T(F)/(\delta N)\).

Now add to \(F\cap S\) every index with phase in
\((1/3+2\delta,2/3-2\delta)\), obtaining \(\overline F\).
In a triple of phases in \(S\) with no wrap, the two inputs
are at most \(1/3+2\delta\) and the output is at least \(2/3-2\delta\).
In a wrapped triple, the real sum is at least \(4/3-\delta\),
so both inputs are at least \(2/3-2\delta\) and the output at most
\(1/3+2\delta\). No triple uses a point of the added open interval.
This proves \eqref{eq:strip} and both count and mass bounds.
\end{proof}

\section{Density from the original triple count}\label{sec:integer-density}

Put \(a_0=1/3\), \(b_0=2/3\), and \(\phi(n)=\{\omega n\}\) in
\([0,1)\). Suppose \(F\subseteq[N]\) contains all indices with
phase in \([a_0+\delta,b_0-\delta]\) and lies in the indices with
phase in \([a_0-\delta,b_0+\delta]\). All deviations from the
middle third lie in its two endpoint strips. We write \(r_F(z)\)
for the number of ordered pairs in \(F^2\) with sum \(z\), for each
\(z\in[N]\). The factor \(\ind_F(z)\) selects the contributions
to \(T(F)\). Let \(E_M\) be the supremum of interval discrepancies
of the first \(M\) phases, with arbitrary endpoint conventions.

Consider first the model set \(F_0=\{n\in[N]:\phi(n)\in I\}\).
Suppose \(\phi(z)=b_0+x\) with \(0<x<1/3\). The two input
phases of a representation from \(F_0\) sum to a number in
\([2/3,4/3)\). The only number in that range congruent to
\(\phi(z)\) modulo one is \(b_0+x\). For \(1\le n<z\)
the two membership conditions are therefore equivalent to
\(\phi(n)\in[a_0,a_0+x]\). The reflection
\(y\mapsto b_0+x-y\) interchanges the ends of this interval.
For example, an output at \(3/4\) pairs input phases in
\([1/3,5/12]\). Thus
\[
 \big|r_{F_0}(z)-(z-1)x\big|\le E_{z-1}.
\]
The window sits at the left endpoint even though the output
is at the right. For a perturbed set \(F\), we widen it to include
the whole left strip. In the lower bound we count the two input
memberships separately, so their common signed change enters twice.

\begin{lemma}[Representations in a prefix]\label{lem:prefix}
Suppose \(0<\delta\le1/24\). For \(z\in[N]\) put \(M=z-1\).
Let \(L_z,R_z\) be the signed sums of
\(\ind_F(n)-\ind_I(\phi(n))\) over \(1\le n<z\) in the left
and right endpoint strips \([a_0-\delta,a_0+\delta]\) and
\([b_0-\delta,b_0+\delta]\), respectively. For \(|x|\le2\delta\),
\begin{align}
 r_F(z)&\ge (Mx+2L_z-3E_M)_+ && \text{if }\phi(z)=b_0+x, \label{eq:repr-right}\\
 r_F(z)&\ge (-Mx+2R_z-3E_M)_+ && \text{if }\phi(z)=a_0+x. \label{eq:repr-left}
\end{align}
\end{lemma}

\begin{proof}
In the first case restrict to indices of the prefix with phase in
\(W=[a_0-3\delta,a_0+x+3\delta]\). This window is preserved by the
phase reflection \(y\mapsto b_0+x-y\) of \(n\mapsto z-n\).
Let \(P\) be the indices in the window. For binary indicators,
\(\ind_F(n)\ind_F(z-n)\ge\ind_F(n)+\ind_F(z-n)-1\).
The resulting count on \(P\) is at least \(2|F\cap P|-|P|\).
The window contains the whole left strip and lies below the right
strip, so the number in \(F\cap P\) is the count in the baseline
interval \([a_0,a_0+x+3\delta]\) plus \(L_z\).
The baseline interval has length \(x+3\delta\), and \(W\)
has length \(x+6\delta\). Taking the baseline count twice and
subtracting \(|P|\) gives the term \(M x\), with error at most
\(3E_M\). Adding \(2L_z\) and using nonnegativity gives
\eqref{eq:repr-right}.

In the second case use \(W=[b_0+x-3\delta,b_0+3\delta]\)
and the reflection \(y\mapsto1+a_0+x-y\).
The baseline is \([b_0+x-3\delta,b_0)\), and the signed change
is \(R_z\). The lengths \(3\delta-x\) and \(6\delta-x\)
give \eqref{eq:repr-left}. The involution preserves ordered
representations including the fixed point \(n=z-n\).
\end{proof}

\begin{lemma}[Density bound for a rotation set]\label{lem:integer-density}
Assume the strip inclusions above hold with \(0<\delta\le1/100\).
Let \(d_0,K_0\ge0\) be such that on each consecutive block of \(M\)
integers in \([N]\), the discrepancy on any circle interval with
arbitrary endpoints is at most \(Md_0+K_0\). Then
\begin{equation}\label{eq:integer-density}
 \frac{|F|}{N}\le\frac13+\sqrt{I_N(F)}
             +1000\left(d_0+\sqrt{\frac{K_0+2}{N}}\right).
\end{equation}
\end{lemma}

\begin{proof}
We sum the prefix bounds by assigning each output integer a cell
of area \(1/N\). The cells put the integers in phase order within
each block, so the bounds near the two endpoints become integrals.

\paragraph{Arranging the outputs.}
Partition \([N]\) into consecutive blocks \(J_j\) of lengths \(m_j\).
Write \(m_*=\min_jm_j\), \(h=\max_j m_j/N\),
and \(D=d_0+(K_0+2)/m_*\). We first work with a partition
for which \(\delta+2D\le1/24\).
In \(J_j=[A_j+1,A_j+m_j]\cap\Z\), sort the indices by phase.
Assign the \(r\)-th index \(n_{j,r}\) to the cell
\([(r-1)/m_j,r/m_j)\) in the phase coordinate.
The ranks of a point strictly below and at or below its phase
give, by the discrepancy bound,
\[
 \frac{r}{m_j}-d_0-\frac{K_0}{m_j}
 \le \phi(n_{j,r})
 \le\frac{r-1}{m_j}+d_0+\frac{K_0}{m_j}.
\]
In particular each assigned phase differs from the actual one
by at most \(D\), including ties.
At time \(s\in(A_j/N,(A_j+m_j)/N]\) and phase \(y\) in the
\(r\)-th cell, set \(z(s,y)=n_{j,r}\) and
\(A(s,y)=\ind_F(z(s,y))\).
The cell of each integer has area \((m_j/N)/m_j=1/N\).
It carries that integer's original representation count, so
\begin{equation}\label{eq:cell-triples}
 I_N(F)=\int_0^1\int_0^1 A(s,y)\,\frac{r_F(z(s,y))}{N}\,dy\,ds .
\end{equation}

Set \(d=\delta+D\). The new profile \(A(s,\cdot)\) contains
the core \([a_0+d,b_0-d]\) and vanishes off \([a_0-d,b_0+d]\).
Let \(\ell(s),r(s)\) be its signed deviations from the middle
third to the left and right of \(1/2\), integrated in phase.
Write \(L(s)=\int_0^s\ell\), \(R(s)=\int_0^s r\).
They satisfy
\begin{equation}\label{eq:profile-mass}
 |\ell(s)|,|r(s)|\le d,\qquad
 L(1)+R(1)=|F|/N-1/3 .
\end{equation}

\paragraph{Prefix comparison.}
On a complete block, changing from the actual phase to the
assigned phase changes the baseline indicator near either
endpoint only on a phase interval of length at most \(2D\).
If the actual and assigned phase cross the dividing point \(1/2\),
both phases are in the full core and have zero signed change.
Thus the error in either normalized signed sum from a complete
block is at most \(2Dm_j/N\).
The incomplete block contributes at most \(2h\), counting
a partial sum and a partial integral. For \(z\) and \(s\) in the
same time block, writing \(M=z-1\), we obtain
\begin{equation}\label{eq:prefix-errors}
 |L_z/N-L(s)|,\ |R_z/N-R(s)|\le 2sD+2h,\qquad
 E_M/N\le sD+h .
\end{equation}
For the last inequality use \(M/N\le s+h\) and
\(hd_0+K_0/N\le hD\le h\).

For a right output phase \(b_0+x\) with \(-2d\le x\le d\),
the assigned index \(z\) has actual phase \(b_0+x_z\) with
\(|x_z-x|\le D\) and \(|x_z|\le2(\delta+2D)\).
Use Lemma~\ref{lem:prefix} at strip width \(\delta+2D\).
The extra parts of the strips have zero signed deviation, so
the same \(L_z,R_z\) apply. Replacing \((M/N)x_z\) by \(sx\)
costs at most \(sD+h\). The replacement of \(2L_z/N\) costs
at most \(4(sD+h)\), and the discrepancy term is at most
\(3(sD+h)\) by \eqref{eq:prefix-errors}.
Taking \(e(s)=16(sD+h)\) to cover these errors,
\eqref{eq:repr-right} gives
\[
 r_F(z)/N\ge(sx+2L(s)-e(s))_+ .
\]
For a left output phase \(a_0-x\) with \(-2d\le x\le d\), the same
replacement in \eqref{eq:repr-left} gives
\(r_F(z)/N\ge(sx+2R(s)-e(s))_+\).

\paragraph{Endpoint integrals.}
Multiply the two bounds by their output indicators \(A(s,b_0+x)\)
and \(A(s,a_0-x)\) and integrate over the corresponding windows.
The windows are disjoint and the integrand in
\eqref{eq:cell-triples} is nonnegative. Thus the two integrals
together are bounded by \(I_N(F)\), including when a window cuts
through a cell.

In each window the output indicator, as a function of \(x\),
is one to the left of \(-d\) and zero to the right of \(d\),
using these values to extend it to the whole real line.
The weight vanishes for \(x<-2d\) by \eqref{eq:profile-mass},
so the extended integral equals the window integral.
For the right output, the indicator has the same signed
mass relative to \(\ind_{\{x<0\}}\) as \(\ind_{\{x<r(s)\}}\).
The latter minimizes the integral against the increasing
weight \(w(x)=(sx+2L(s)-e(s))_+\).
The difference of the two indicators has zero integral, and
its product with \(w(x)-w(r(s))\) is nonnegative everywhere.
Consequently the right integral at \(s\) is at least
\(s(v_R(s))_+^2/2\), where
\[
 v_R(s)=r(s)+\frac{2L(s)-e(s)}s,\qquad
 v_L(s)=\ell(s)+\frac{2R(s)-e(s)}s .
\]
The left integral is at least \(s(v_L(s))_+^2/2\) by the same
comparison with the reflected coordinate \(a_0-x\).

The two endpoint terms contain the accumulated changes
\(L(s),R(s)\) and their derivatives \(\ell(s),r(s)\).
On adding the terms, the accumulated changes occur with
coefficient \(2/s\). Multiplication by \(s^2\) therefore gives
\[
 s^2(v_R(s)+v_L(s))=(s^2(L(s)+R(s)))'-2s e(s).
\]
The term at zero vanishes since \(L(s),R(s)=O(s)\). Integrate, replace both
terms on the left by their positive parts, and use
Cauchy--Schwarz with
\(\tfrac12\int_0^1s((v_L)_+^2+(v_R)_+^2)\,ds\le I_N(F)\).
The second factor in Cauchy--Schwarz is
\((2\int_0^1s^3\,ds)^{1/2}=1/\sqrt2\), giving
\[
 L(1)+R(1)\le\sqrt{I_N(F)}+2\int_0^1s e(s)\,ds
            =\sqrt{I_N(F)}+(32/3)D+16h.
\]

\paragraph{Block length.}
Write \(\sigma=d_0+\sqrt{(K_0+2)/N}\).
If \(\sigma\ge1/1000\), \eqref{eq:integer-density} is immediate.
Otherwise take blocks of length \(M=\lceil\sqrt{(K_0+2)N}\rceil\),
merging the final short block with its predecessor. We have
\(M<N/2\), \(m_*\ge M\), and \(h\le2M/N\le2\sigma+\sigma^2\).
Also \(D\le\sigma\), so \(\delta+2D\le1/24\) and
\((32/3)D+16h\le1000\sigma\). With \eqref{eq:profile-mass}
this proves the claimed bound.
\end{proof}

\section{Proof of the small-shift estimate}\label{sec:proof-shifts}

We prove a stronger shift assertion first. If
\begin{equation}\label{eq:dense-F}
 |F|\ge (1/3+7\eta/8)N,\qquad I_N(F)\le a\eta^2 ,
\end{equation}
then \(F\cap(F+t)\ne\varnothing\) for each positive multiple \(t\)
of \(Q_\eta\) with \(t\le\kappa_\eta N\).
The model accuracy is fixed once. At the integer counting stage
the growing cutoff \(H_\eta\), the restriction \(t/N\le K\eta^2\),
and the triple count in \eqref{eq:dense-F} will each give an error
on the scale of \(\eta\).

\paragraph{Constants.}
Set \(\delta_0=1/1000\), \(\gamma=\delta_0/100\),
\(\epsilon_0=\gamma^2/4096\), and \(b_*=1/4\).
Take \(D,J_*\) from Proposition~\ref{prop:structure} at
accuracy \(\epsilon_0/4\).
Choose \(\xi=\min(1/2,\sqrt{J_*}/2)\) and take \(L\)
from Lemma~\ref{lem:return} at this \(\xi\).
Fix
\[
 \begin{aligned}
 H_0 &= 32000,\\
 K &= \min\{1/(272L),\, b_* /(256\cdot2000^2),\,1/16\},\\
 a &= \min\{1/256,\,\delta_0/128,\,\gamma^2/1024,\,J_*/12\},\\
 \eta_0 &= \min\{1/12,\,H_0/D,\,\epsilon_0,\,\gamma/1024\}.
 \end{aligned}
\]
In particular \(H_\eta\ge D\) and \(\kappa_\eta\le\eta/16\)
for \(0<\eta\le\eta_0\). The constants precede the choice
of \(\eta,N,F,t\).

\paragraph{A coarse model for \(F\).}
Suppose \(F\) in \eqref{eq:dense-F} is independent under an
eligible \(t\). Let \(b=\lceil16/\eta\rceil\le17/\eta\).
Then \(btL/N\le17KL\eta\le\eta/16\).
Lemma~\ref{lem:return} gives a multiple \(q\) of \(bt\) with
\(q/N,t/q\le\eta/16\) and \(T_q(F)\le2T(F)+\xi^2N^2\).
By Lemma~\ref{lem:sampler} the resulting profile on \(G=\Z/q\Z\)
is independent under \(v=t\bmod q\), with
\[
 d(h)\ge1/3+3\eta/4,\qquad J(h)\le3a\eta^2+\xi^2\le J_*.
\]
Since \(Q(D)\mid Q_\eta\mid t\), we have \(v\in Q(D)G\).
Proposition~\ref{prop:structure} gives \(\chi,\lambda\) at
error at most \(\epsilon_0/4\), with \(\|\chi(v)\|_\T\ge5/16\).
Lemma~\ref{lem:unfold} yields a rotation of frequency \(\omega\)
and offset \(\beta\) with approximation error at most
\(\epsilon_0/(4(1-t/q))+q/N\le\epsilon_0\). It also gives
\[
 \|\omega t\|_\T>5/16-\frac{4\eta}{16-\eta}>b_* .
\]

\paragraph{Counting on the integers.}
Lemma~\ref{lem:discrepancy} gives discrepancy at most
\(\mathcal D=N/(H_\eta+1)+3t/(2b_*)\) on every block of \([N]\).
Thus \(\mathcal D/N\le\eta/H_0+6K\eta^2\le3\eta/2\le\gamma/64\).
The existence of \(t\ge1\) with \(t/N\le K\eta^2\) ensures \(N\ge16\).
Apply Lemma~\ref{lem:offset} at \(\delta=\delta_0\) to obtain
\(\overline F\) having the required strip inclusions at width
\(3\delta_0\le1/100\) and \(T(\overline F)\le T(F)\).
For the block discrepancy in Lemma~\ref{lem:integer-density},
take \(d_0=1/(H_\eta+1)\) and \(K_0=3t/(2b_*)\).
Since \(K_0+2\le5t/(2b_*)\), the density and deletion bounds give
\[
 \frac{|F|}{N}\le\frac13+\sqrt a\,\eta
        +\frac{2000}{H_\eta+1}
        +2000\sqrt{K/b_*}\,\eta
        +\frac{8a\eta^2}{\delta_0}.
\]
Each of the four errors is at most \(\eta/16\) by the choices
of constants. This contradicts \eqref{eq:dense-F}.

For \(A\) in the proposition, suppose \(L_A(t)\le\eta N/8\)
at some eligible \(t\). Delete all \(x\in A\) for which
\(x+t\in A\). The remaining set is independent under \(t\),
has the density in \eqref{eq:dense-F}, and its count can only
decrease. The contradiction proves \(L_A(t)>\eta N/8\) for
every eligible \(t\). Sum over those \(t\) belonging to \(A\)
in \eqref{eq:shift-triples} to obtain \eqref{eq:origin}.

Finally \(Q(2m)/Q(m)\) divides \(\binom{2m}{m}\).
For each prime whose exponent increases, exactly one new prime
power lies in \((m,2m]\). It contributes one to the valuation of
the binomial coefficient, with all other summands nonnegative.
Therefore \(\log Q(2m)-\log Q(m)\le2m\log2\).
Summing at powers of two and bounding \(H_\eta\) by the next
power of two gives \(\log Q_\eta=O(H_\eta)=O(1/\eta)\).
This completes the proof of Proposition~\ref{prop:shifts} and
Theorem~\ref{thm:main}.

\appendix
\section{The slice-density inequality}\label{app:slice-density}

\begin{proof}[Proof of Lemma~\ref{lem:triangle}]
Write \(E(u,w)=g(u)+g(w)+g(u+w)\).
Put \(K=[1/4,2/3]\), \(J=[1/4,1/3]\), \(B=[1/16,1/8]\).
For \(u\in K\), \(w\in J\), \(y\in B\), cancellation gives
\[
 2(g(u)-g(w))=E(u,w-y)+E(u-y,y)-E(u-y,w)-E(y,w-y).
\]
All four sums are in \([1/4,1]\). Integrate the absolute values
over \(K\times J\times B\), using the displayed input variables
of each residual as integration coordinates. The changes of
variables have absolute Jacobian one, and the remaining fibers
have lengths at most \(|B|,|J|,|B|,|K|\), respectively. Since the
weight in each residual is at least \(3/4\), this gives
\[
 2|B|\int_K\int_J |g(u)-g(w)|\,dw\,du
 \le \tfrac43(2|B|+|J|+|K|)R=\tfrac56 R .
\]
With \(c=|J|^{-1}\int_Jg\), averaging over \(w\) yields
\(\int_K|g-c|\le80R\). On \(J^2\), all three times lie in \(K\)
and the weight is at least \(3/2\). Each of the three deviations
from \(c\) integrates to at most \(|J|\int_K|g-c|\), including
the sum variable, whose unnormalized density is at most \(|J|\).
Hence \(3|c||J|^2\le 2R/3+3|J|\cdot80R\),
so \(|c|\le992R\) and \(G_K=\int_K|g|\le (1480/3)R\).

For \(0<t<1/4\), average
\(|g(t)|\le |E(t,w)|+|g(w)|+|g(t+w)|\)
over \(w\in[1/3,5/12]\). The residual weight is at least one,
both other times are in \(K\), and their integration fibers
have lengths at most \(1/4\) and \(1/12\). Dividing by \(1/12\)
gives \(\int_0^{1/4}|g|\le12R+4G_K\).
For \(2/3<t<1\), average
\(|g(t)|\le |E(u,t-u)|+|g(u)|+|g(t-u)|\)
over \(|u-t/2|\le1/24\).
Both inputs are in \(K\), the residual weight is at least two,
and the \(t\)-fiber for each fixed input has length at most \(1/6\).
This gives \(\int_{2/3}^1|g|\le6R+4G_K\).
The three integrals sum to at most \(18R+9G_K\le4500R\).
\end{proof}

\subsection*{Acknowledgments}
We used GPT-6 Astra through Codex to explore proof strategies and to
develop and check arguments. Codex also assisted in preparing and
revising the manuscript. We used
\href{https://github.com/zxl0714/proof-by-gacha}{Gacha}, a tool we
developed, to manage research workspaces and runs.

\end{document}